\documentclass[12pt]{amsart}

\usepackage{amsmath,xspace,amssymb,mathrsfs}
\usepackage{color}

\input xy
\xyoption{all}
\xyoption{2cell}
\UseAllTwocells
\CompileMatrices

\newcommand{\Spec}{\operatorname{Spec}}
\renewcommand{\phi}{\varphi}

\newcommand{\CR}{\operatorname{C-Ring}}

\newcommand{\D}{\operatorname{D}}
\newcommand{\Clop}{\operatorname{Clop}}
\newcommand{\Fin}{\operatorname{Fin}}
\newcommand{\Char}{\operatorname{Char}}
\newcommand{\Top}{\operatorname{Top}}

\newcommand{\Set}{\operatorname{Set}}
\newcommand{\Ring}{\operatorname{Ring}}
\newcommand{\Mor}{\operatorname{Mor}}

\newcommand{\Sp}{\operatorname{Sp}}
\newcommand{\V}{\operatorname{V}}

\newtheorem{proposition}{Proposition}[section]
\newtheorem{lemma}[proposition]{Lemma}

\newtheorem{corollary}[proposition]{Corollary}
\newtheorem{theorem}[proposition]{Theorem}

\theoremstyle{definition}

\newtheorem{remark}[proposition]{Remark}

\begin{document}

\title[Stone type representations and dualities]{Stone type representations and dualities by power set ring}

% On the structure and application of power set ring

\author[A. Tarizadeh and Z. Taheri]{Abolfazl Tarizadeh and Zahra Taheri}
\address{Department of Mathematics, Faculty of Basic Sciences, University of Maragheh \\
P. O. Box 55136-553, Maragheh, Iran.
 }
\email{ebulfez1978@gmail.com, zta.pardafan@gmail.com}

\date{}

\begin{abstract} In this paper, it is shown that the Boolean ring of a commutative ring is isomorphic to the ring of clopens of its prime spectrum. In particular, Stone's Representation Theorem is generalized. The prime spectrum of the Boolean ring of a given ring $R$ is identified with the Pierce spectrum of $R$. The discreteness of prime spectra is characterized. It is also proved that the space of connected components of a compact space $X$ is isomorphic to the prime spectrum of the ring of clopens of $X$. As another major result, it is shown that a morphism of rings between complete Boolean rings preserves suprema if and only if the induced map between the corresponding prime spectra is an open map.
\end{abstract}

\maketitle

\tableofcontents

\section{Introduction}

In 1936, Marshall Stone published a long paper \cite{Stone} that whose main result was that every Boolean ring is isomorphic to a certain subring of a power set ring. This Stone Representation Theorem has been the starting point for a whole new field of study, nowadays called Stone duality. It is worth mentioning that ``Boolean ring'' and ``Boolean algebra'' are the equivalent structures, i.e., every Boolean ring can be made canonically into a Boolean algebra (and vice versa). This is a well known fact, for its proof and further studies on these topics we refer the interested reader to the References, especially to \cite{Borceux-Janelidze}, \cite{Halmos 2}, \cite{Halmos} and \cite{Johnstone}. \\
In this paper we prove that if $R$ is a commutative ring, then we have the canonical isomorphism of rings:
$$\mathcal{B}(R)\simeq\Clop\big(\Spec(R)\big)$$
here $\mathcal{B}(R)$ denotes the ring of idempotents of $R$; we call $\mathcal{B}(R)$ the Boolean ring of $R$. This is Theorem \ref{Theorem I}. Note that the Stone space of a Boolean algebra (see e.g. \cite[Proposition 4.1.11]{Borceux-Janelidze}) and the prime spectrum of the corresponding Boolean ring are canonically homeomorphic. In fact, this homeomorphism takes each ultrafilter $F$ in a given Boolean algebra $R$ into $R\setminus F$ which is a maximal ideal of the corresponding Boolean ring. By taking into account this identification, then the original Stone's Representation Theorem is a especial case of the above isomorphism, see Corollary \ref{Corollary III}. \\
Indeed, considering $\Clop(X)$ as a ring (where $X$ is a topological space) is the simple idea leading to some results of this paper. In Theorem \ref{Theorem V}, we prove the topological version of  Theorem \ref{Theorem I} which states that if $X$ is a compact space then we have the canonical isomorphism of topological spaces: $$\pi_{0}(X)\simeq\Spec\big(\Clop(X)\big)$$ where $\pi_{0}(X)$ denotes the space of connected components of $X$ which carries the quotient topology. If moreover $X$ is totally disconnected, then we easily get the following well known isomorphism which can be considered it as the topological version of the Stone's Representation Theorem: $$X\simeq\Spec\big(\Clop(X)\big).$$
In Sections 3 and 4 we prove several structural results on the Boolean ring of a given ring and its prime spectrum. Especially, the structure of prime ideals of the Boolean ring of a given ring $R$ is identified in terms of the connected components of $\Spec(R)$. In fact, in Theorem \ref{Theorem PSp}, we obtain the following canonical isomorphism of topological spaces: $$\Spec\mathcal{B}(R)\simeq\pi_{0}\big(\Spec(R)\big).$$
Then as an application, for a given ring $R$, it is shown that $X=\Spec(R)$ is discrete if and only if $\mathcal{B}(R)\simeq\mathcal{P}(X)$, see Theorem \ref{Thm disc. char}. \\
Sections 5 and 6 investigate the structure of Boolean rings, especially complete Boolean rings. Theorem \ref{Theorem DII} gives a canonical bijection from the set of ring maps between two Boolean rings onto the set of continuous maps between their corresponding prime spectra. \\
In Theorem \ref{Theorem Guram}, we prove a technical result which states that a morphism of rings between complete Boolean rings preserves suprema if and only if the induced map between the corresponding prime spectra is an open map. This observation leads us to various non trivial results. Especially among them, we obtain a Stone type duality which states that the category of complete Boolean rings is the dual of the category of compact extremally disconnected spaces, see Theorem \ref{Thm anti-equivalence}. \\
In \S7, we have made further progresses in the understanding the structure of power set ring. Also some contributions to the field of fixed-point theorems are provided.

\section{Preliminaries}

We collect in this section some basic background for the reader's convenience. \\
In this paper, all rings are commutative. If $f$ is an element of a ring $R$, then $\D(f)=\{\mathfrak{p}\in\Spec(R): f\notin\mathfrak{p}\}$ and $\V(f)=\Spec(R)\setminus\D(f)$. \\
If $X$ is a set then its power set $\mathcal{P}(X)$ together with the symmetric difference $A+B=(A\cup B)\setminus (A\cap B)$ as the addition and the intersection $A.B=A\cap B$ as the multiplication forms a commutative ring whose zero and unit are respectively the empty set and the whole set $X$. The ring $\mathcal{P}(X)$ is called the \emph{power set ring} of $X$. If $\phi:X\rightarrow Y$ is a function then the map $\mathcal{P}(\phi):\mathcal{P}(Y)\rightarrow\mathcal{P}(X)$ defined by $A\rightsquigarrow\phi^{-1}(A)$ is a morphism of rings. In fact, the assignments $X\rightsquigarrow\mathcal{P}(X)$ and $\phi\rightsquigarrow\mathcal{P}(\phi)$ form a faithful contravariant functor from the category of sets to the category of commutative rings. We call it the \emph{power set functor}. \\
Remember that a ring is called a Boolean ring if each element is an idempotent. Power set ring $\mathcal{P}(X)$ is a typical example of Boolean rings. It is easy to see that in a Boolean ring, every prime ideal is a maximal ideal. \\
By a compact space we mean a quasi-compact and Hausdorff topological space. The number of elements (cardinal) of a set $X$ is denoted by $|X|$. \\
If $\phi:R\rightarrow R'$ is a morphism of rings then the induced map $\Spec(R')\rightarrow\Spec(R)$ defined as $\mathfrak{p}\rightsquigarrow\phi^{-1}(\mathfrak{p})$
is denoted by $\Spec(\phi)$. \\
If $A$ and $B$ are objects of a category $\mathscr{C}$, then as usual by $\Mor_{\mathscr{C}}(A,B)$ we mean the set of morphisms of $\mathscr{C}$ from $A$ to $B$. \\
By a \emph{join-complete poset} we mean a poset $A$ such that the supremum of every subset $(x_{k})\subseteq A$ exists in A, this element is denoted by $\bigvee\limits_{k}x_{k}$. For example the posets $\mathcal{P}(X)$ and $\Spec(\mathbb{Z})$ ordered by the inclusion, are join-complete and not join-complete, respectively. \\
Dually, by a \emph{meet-complete poset} is meant a poset $A$ such that the infimum of every subset $(x_{k})\subseteq A$ exists in A, this element is denoted by $\bigwedge\limits_{k}x_{k}$. Recall that a poset which is both join-complete and meet-complete is called a \emph{complete lattice}. \\
It is easy to see that an order preserving map $f:A\rightarrow B$ between the join-complete posets preserves suprema if and only if there exists an order preserving map $g:B\rightarrow A$ such that the pair $(f,g)$ is a Galois connection, i.e., $f(a)\leqslant b\Leftrightarrow a\leqslant g(b)$ for all $a\in A$ and $b\in B$. In this case, such function $g$ is unique. In fact, if $f$ preserves suprema, then $g:B\rightarrow A$ is defined as $g(b)=\bigvee\limits_{f(a)\leqslant b}a$. \\
Dually, an order preserving map $f:A\rightarrow B$ between the meet-complete posets preserves infima if and only if there exists an order preserving map $g:B\rightarrow A$ such that the pair $(g,f)$ is a Galois connection. If $f$ preserves infima, then $g(b)=\bigwedge\limits_{f(a)\geqslant b}a$.

\section{Generalization of Stone Representation Theorem}

In this section, Stone's Representation Theorem is generalized from Boolean rings to arbitrary commutative rings. \\
If $X$ is a topological space then by $\Clop(X)$ we mean the set of clopen (both open and closed) subsets of $X$. It has the canonical ring structure. In fact, $\Clop(X)$ is a subring of $\mathcal{P}(X)$. If $\phi:X\rightarrow Y$ is a continuous map of topological spaces then the map $\Clop(\phi):\Clop(Y)\rightarrow\Clop(X)$ given by $A\rightsquigarrow \phi^{-1}(A)$ is a morphism of rings.
In fact, the assignments $X\rightsquigarrow\Clop(X)$ and $\phi\rightsquigarrow\Clop(\phi)$ form a contravariant functor from the category of topological spaces to the category of Boolean rings. We call it the \emph{Clopen functor}. \\
Let $R$ be a ring and $\mathcal{B}(R)=\{e\in R: e=e^{2}\}$. Then it is easy to see that the set $\mathcal{B}(R)$ by the new operation $\oplus$ defined by $e\oplus e':=e+e'-2ee'$ as the addition admits a ring structure whose multiplication is the multiplication of $R$. Note that $\mathcal{B}(R)$ is not necessarily a subring of $R$. In fact, $\mathcal{B}(R)$ is a subring of a non-zero ring $R$ if and only if $\Char(R)=2$. We call $\mathcal{B}(R)$ the \emph{Boolean ring} (or, \emph{Booleanization}) of $R$. Clearly $\mathcal{B}(R)=R$ if and only if $R$ is Boolean. We have then:

\begin{theorem}$($Generalized Stone Representation Theorem$)$\label{Theorem I} Let $R$ be a ring and $X=\Spec(R)$. Then the map $\mathcal{B}(R)\rightarrow\Clop(X)$ given by $e\rightsquigarrow\D(e)$ is an isomorphism of rings.
\end{theorem}

{\bf Proof.} It is well known that the above map is bijective, see e.g. \cite[Proposition 3.1]{Tarizadeh}. Clearly $\D(ee')=\D(e)\cap\D(e')$ and $\D(1)=\Spec(R)$. It remains to show that $\D(e\oplus e')=\D(e)+\D(e')$. If $\mathfrak{p}\in\D(e\oplus e')$ then clearly $\mathfrak{p}\in\D(e)\cup\D(e')$. But $1-(e\oplus e')\in\mathfrak{p}$, because $e\oplus e'$ is an idempotent. So $ee'=ee'\big(1-(e\oplus e')\big)\in\mathfrak{p}$. This yields that $\mathfrak{p}\in\D(e)+\D(e')$. Conversely, assume that $\mathfrak{p}\in\big(\D(e)\cup\D(e')\big)\setminus\D(ee')$. If $e\oplus e'\in\mathfrak{p}$ then $e+e'-ee'\in\mathfrak{p}$. It follows that $e=e(e+e'-ee')\in\mathfrak{p}$ and $e'=e'(e+e'-ee')\in\mathfrak{p}$. But this is a contradiction. Hence, $\mathfrak{p}\in\D(e\oplus e')$. $\Box$

\begin{corollary}$($Stone Representation Theorem$)$\label{Corollary III} If $R$ is a Boolean ring then the map $R\rightarrow\Clop\big(\Spec(R)\big)$ given by $e\rightsquigarrow\D(e)$ is an isomorphism of rings.
\end{corollary}

{\bf Proof.} It follows from Theorem \ref{Theorem I}. $\Box$

The above result, in particular, tells us that a Boolean ring $R$ is finite if and only if $\Spec(R)$ is finite.

\begin{corollary}\label{Corollary ix injective b} Let $R$ be a ring and $X=\Spec(R)$. Then the canonical map $\mathcal{B}(R)\rightarrow\mathcal{P}(X)$ given by $e\rightsquigarrow \D(e)$ is an injective morphism of rings.
\end{corollary}

{\bf Proof.} It is an immediate consequence of Theorem \ref{Theorem I}. $\Box$

\begin{remark} Let $R$ be a ring. By \cite[Theorem 6]{Hochster} or \cite[Theorem 3.20]{Tarizadeh 2}, there exists a ring $S$ and a homeomorphism $\phi:X_{\mathcal{Z}}=\big(\Spec(R),\mathcal{Z}\big)\rightarrow
Y_{\mathcal{F}}=\big(\Spec(S),\mathcal{F}\big)$ such that if $\mathfrak{p}\subseteq\mathfrak{q}$ are prime ideals of $R$, then $\phi(\mathfrak{q})\subseteq\phi(\mathfrak{p})$ where $\mathcal{Z}$ and $\mathcal{F}$ denote the Zariski and flat topologies, respectively. We call $S$ a \emph{prime-inverse ring} of $R$. Note that prime-inverse ring is not unique. In general, there is no way to find a ring map between $R$ and $S$. In spite of this, we have the following result.
\end{remark}

\begin{corollary} Let $R$ be a ring. If $S$ is a prime-inverse ring of $R$, then the Boolean rings $\mathcal{B}(R)$ and $\mathcal{B}(S)$ are isomorphic.
\end{corollary}

{\bf Proof.} Clearly $\Clop(\phi):\Clop(Y_{\mathcal{F}})\rightarrow
\Clop(X_{\mathcal{Z}})$ is an isomorphism of rings, for $\phi$ see the above remark. By \cite[Corollary 3.12]{Tarizadeh 2}, $\Clop(Y_{\mathcal{F}})=\{\D(e): e\in\mathcal{B}(S)\}=\Clop(Y_{\mathcal{Z}})$.
Then apply Theorem \ref{Theorem I}. $\Box$

If $\phi:R\rightarrow R'$ is a morphism of rings then the map $\mathcal{B}(\phi):\mathcal{B}(R)\rightarrow\mathcal{B}(R')$ given by $e\rightsquigarrow\phi(e)$ is also a morphism of rings. In fact,
the assignments $R\rightsquigarrow\mathcal{B}(R)$ and $\phi\rightsquigarrow\mathcal{B}(\phi)$ define a covariant functor from the category of commutative rings to the category of Boolean rings. We call it the Booleanization functor. In fact, the isomorphism of Theorem \ref{Theorem I} is functorial. Namely, the Booleanization functor is isomorphic to $\Clop\circ\Spec$ where the contravariant functor $\Spec:\CR\rightarrow\Top$ is the usual spectrum functor ($\CR$ denotes the category of commutative rings).

\section{Booleanization versus Pierce spectrum}

Remember from \cite[\S 2]{Tarizadeh-Aghajani} or \cite[\S3]{Tarizadeh 2} that an ideal of a ring $R$ is called a regular ideal if it is generated by a set of idempotents of $R$. If $I$ is a regular ideal of a ring $R$ and $f\in I$, then $f=fe$ for some idempotent $e\in I$, see \cite[Lemma 8.4]{Tarizadeh-Aghajani}. \\
Each maximal element of the set of proper and regular ideals of $R$ is called a max-regular ideal of $R$. The set of max-regular ideals of $R$ is called the \emph{Pierce spectrum} of $R$ and denoted by $\Sp(R)$. It is a compact totally disconnected space whose base opens are precisely of the form $U_{e}=\{M\in\Sp(R): e\notin M\}$ where $e\in R$ is an idempotent. It can be seen that the map $\Sp(R)\rightarrow\pi_{0}\big(\Spec(R)\big)$ given by $M\rightsquigarrow\V(M)$ is a homeomorphism.\\
As stated in \S3, the ring $\mathcal{B}(R)$ is not necessarily a subring of $R$. In spite of this, we have the following topological identification.

\begin{theorem}\label{Theorem PSp} Let $R$ be a ring. Then $\Spec\mathcal{B}(R)$ is canonically homeomorphic to the Pierce spectrum $\Sp(R)$.
\end{theorem}

{\bf Proof.} We show that the map $\mu:\Spec\mathcal{B}(R)\rightarrow\Sp(R)$ given by $P\rightsquigarrow(P)$ is a homeomorphism where $(P)$ is an ideal of $R$ generated by the subset $P$. If $P$ is a prime ideal of $\mathcal{B}(R)$ then $(P)\neq R$, because if $1\in(P)$ then we may write $1=\sum\limits_{i=1}^{m}r_{i}e_{i}$ where $r_{i}\in R$ and $e_{i}\in P$ for all $i$. Setting $e'=\prod\limits_{k=1}^{m}(1-e_{k})$. We have $e'=1e'=\sum\limits_{i=1}^{m}r_{i}e_{i}e'=0$. But $1-e_{k}\in\mathcal{B}(R)$ for all $k$. Thus $1-e_{k}\in P$ for some $k$. It follows that $1=(1-e_{k})+e_{k}\in P$ which is a contradiction. In fact, it is easy to see that $(P)$ is a max-regular ideal of $R$. By a similar argument as applied in the above, it is shown that the map $\mu$ is injective. It is also continuous, since $\mu^{-1}(U_{e})=\D(e)=\{P\in\Spec\mathcal{B}(R): e\notin P\}$ for all $e\in\mathcal{B}(R)$. Thus $\mu$ is a closed map, because $\Spec\mathcal{B}(R)$ is quasi-compact and $\Sp(R)$ is Hausdorff. It remains to show that this map is surjective. Let $M$ be a max-regular ideal of $R$. The ideal $P=(e\in M: e=e^{2})$ of $\mathcal{B}(R)$ generated by the subset $M\cap\mathcal{B}(R)$ is a proper ideal. In fact, the set $P$ is contained in $M$, because if $e\in P$ then we may write $e=e_{1}e'_{1}\oplus\cdots\oplus e_{n}e'_{n}\in Re'_{1}+\cdots+Re'_{n}\subseteq M$ where $e_{k}\in\mathcal{B}(R)$ and $e'_{k}\in M\cap\mathcal{B}(R)$ for all $k$. Assume there exist $e,e'\in\mathcal{B}(R)$ such that $ee'\in P$. Then $ee'\in M$. Suppose $e\notin M$ and $e'\notin M$. This yields that $M+Re=R$ and $M+Re'=R$, since $M+Re$ and $M+Re'$ are regular ideals. Thus $1=m+re$ and $1=m'+r'e'$ where $m,m'\in M$. This implies that $1=mm'+mr'e'+rem'+rr'ee'\in M$, a contradiction. Thus we may assume that, say, $e\in M$. So $e\in P$. Hence, $P$ is a prime ideal of $\mathcal{B}(R)$. In the above we observed that $(P)$ is a max-regular ideal of $R$. Clearly $(P)\subseteq M$. Therefore $\mu(P)=(P)=M$. $\Box$

The above theorem, in particular, yields that a ring $R$ is zero dimensional if and only if $\Spec(R)\simeq\Spec\mathcal{B}(R)$. \\
If $R$ is a finite ring or more generally an Artinian ring, then $\Spec(R)$ is discrete. In fact, for a given ring $R$, then $\Spec(R)$ is discrete if and only if $R$ is zero dimensional with finitely many primes.
The following result gives us further characterizations of the discreteness of prime spectra.

\begin{theorem}\label{Thm disc. char} Let $R$ be a ring and $X=\Spec(R)$. Then the following statements are equivalent. \\
$\mathbf{(i)}$ The Zariski topology of $X$ is discrete. \\
$\mathbf{(ii)}$ The canonical ring map $\mathcal{B}(R)\rightarrow\mathcal{P}(X)$ given by $e\rightsquigarrow \D(e)$ is surjective. \\
$\mathbf{(iii)}$ The ring $\mathcal{B}(R)$ is isomorphic to $\mathcal{P}(X)$.
\end{theorem}

{\bf Proof.} $\mathbf{(i)}\Rightarrow\mathbf{(ii)}:$ If $X$ is discrete then $\Clop(X)=\mathcal{P}(X)$. Thus by Theorem \ref{Theorem I}, the above map is surjective. \\
$\mathbf{(ii)}\Rightarrow\mathbf{(iii)}:$ See Corollary \ref{Corollary ix injective b}. \\
$\mathbf{(iii)}\Rightarrow\mathbf{(i)}:$ It suffices to show that $X$ is a finite set. Because in this case, $\mathcal{P}(X)$ will also be a finite set and the correspondences $\Clop(X)\simeq\mathcal{B}(R)\simeq\mathcal{P}(X)$ would imply that $\Clop(X)=\mathcal{P}(X)$, since a finite set has no proper subset with the same cardinality. Suppose $X$ is an infinite set with the cardinality $\kappa$. Using Theorem \ref{Theorem PSp}, then we have $\Spec\mathcal{P}(X)\simeq\Spec\mathcal{B}(R)\simeq\Sp(R)$.
It is well known that $\Spec\mathcal{P}(X)$ is the Stone-\v{C}ech compactification of the discrete space $X$. So by \cite[Theorem 3.58]{Hindman-Strauss} or by \cite[Theorem on p.71]{Walker} we have
$|\Sp(R)|=2^{2^{\kappa}}$. On the other hand, $\Sp(R)$ is homeomorphic to the space of connected components of $X$ in the Zariski topology. Therefore $|\Sp(R)|\leqslant\kappa$. But this is in contradiction with Cantor's theorem (i.e., $\kappa<2^{\kappa}$ for all cardinals $\kappa$). Hence, $X$ is finite. $\Box$

Let $I$ be an ideal of a ring $R$. Then $I^{\ast}=(e\in I: e=e^{2})$ is the largest regular ideal of $R$ which is contained in $I$. We call $I^{\ast}$ the \emph{regular part} of $I$. If $\phi:R\rightarrow R'$ is a morphism of rings and $M$ a max-regular ideal of $R'$, then $\big(\phi^{-1}(M)\big)^{\ast}$ is a max-regular ideal of $R$. The map $\Sp(\phi):\Sp(R')\rightarrow\Sp(R)$ given by $M\rightsquigarrow\big(\phi^{-1}(M)\big)^{\ast}$ is continuous, because $\Sp(\phi)^{-1}(U_{e})=U_{\phi(e)}$ for all $e\in\mathcal{B}(R)$.
If $\psi:R'\rightarrow R''$ is a second morphism of rings, then it is not hard to see that $\Sp(\psi\circ\phi)=\Sp(\phi)\circ\Sp(\psi)$. In fact, the assignments $R\rightsquigarrow\Sp(R)$ and $\phi\rightsquigarrow\Sp(\phi)$ form a contravariant functor from the category of commutative rings to the category of compact totally disconnected spaces. We call it the \emph{Pierce functor}. Then we have the following result.

\begin{proposition} The Pierce functor is isomorphic to the contravariant functor $\Spec\circ\mathcal{B}$.
\end{proposition}

{\bf Proof.} It is easy to see that if $\phi:R\rightarrow R'$ is a morphism of rings then the following diagram is commutative: \\ $$\xymatrix{\Spec\mathcal{B}(R')
\ar[r]^{\:\:\:\:\:\:\:\:\:\mu_{R'}}\ar[d]^{F(\phi)}&\Sp(R')\ar[d]^{\Sp(\phi)}\\
\Spec\mathcal{B}(R)\ar[r]^{\:\:\:\:\:\:\:\:\:\mu_{R}}&\Sp(R)}$$ where $\mu_{R}$ and $\mu_{R'}$ are the canonical isomorphisms (see Theorem \ref{Theorem PSp}) and $F=\Spec\circ\mathcal{B}$. $\Box$

\section{Power set ring and Stone duality}

If $A\in R=\mathcal{P}(X)$ then $\mathcal{P}(A)$ is an ideal of $R$ and the quotient ring $R/\mathcal{P}(A)$ is canonically isomorphic to the ring $\mathcal{P}(A^{c})$ where $A^{c}=X\setminus A$. Especially, if $x\in X$ then $\mathfrak{m}_{x}:=\mathcal{P}(X\setminus\{x\})$
is a maximal ideal of $R$, since $R/\mathfrak{m}_{x}\simeq\mathbb{Z}_{2}=\{0,1\}$. \\
It is easy to see that in a ring $R$, if an ideal is generated by a finite set of idempotents then it is generated by one idempotent. In fact, if $e,e'\in R$ are idempotents then the ideal $(e,e')$ is generated by the idempotent $e+e'-ee'$. In particular, we have the following result.

\begin{proposition}\label{Lemma DI} If $A_{1},\ldots,A_{n}\in\mathcal{P}(X)$ then $(A_{1},\ldots,A_{n})=\mathcal{P}(\bigcup\limits_{k=1}^{n}A_{k})$.
\end{proposition}

{\bf Proof.} First note that for each $A\in\mathcal{P}(X)$ then $(A)=\mathcal{P}(A)$. Because clearly $(A)\subseteq\mathcal{P}(A)$. Conversely, if $C\in\mathcal{P}(A)$ then $C=C\cap A=CA\in (A)$. Now if $B\in\mathcal{P}(X)$ then the ideal $(A,B)$ is generated by the element $A+B-AB=A\cup B$. Thus $(A,B)=\mathcal{P}(A\cup B)$.  $\Box$

\begin{remark}\label{Remark ii subring new} Let $R$ be a subring of $\mathcal{P}(X)$ for some set $X$. Also let $I$ be an ideal of $R$ and $A,B\in I$. Then $A\cup B=A+B+A\cap B\in I$.
\end{remark}

The following result can be viewed as the topological version of Theorem
\ref{Theorem I}.

\begin{theorem}\label{Theorem V} Let $X$ be a compact space and $R=\Clop(X)$. Then the map $\phi:\pi_{0}(X)\rightarrow\Spec(R)$ given by $[x]\rightsquigarrow\mathfrak{m}_{x}\cap R$ is a homeomorphism.
\end{theorem}

{\bf Proof.} We have $[x]=[y]$ if and only if $\mathfrak{m}_{x}\cap R=\mathfrak{m}_{y}\cap R$, because it is well known that if $X$ is a compact space then for each $x\in X$, the connected component $[x]$ is the intersection of all clopens of $X$ which contain $x$. Thus $\phi$ is well-defined and injective. If $A\in R$ then $\phi^{-1}\big(\D(A)\big)$ is an open of $\pi_{0}(X)$, since $\pi^{-1}\Big(\phi^{-1}\big(\D(A)\big)\Big)=A$ where $\pi:X\rightarrow\pi_{0}(X)$ is the canonical projection. Hence, $\phi$ is continuous. It is also a closed map, because $\pi_{0}(X)$ is quasi-compact and $\Spec(R)$ is Hausdorff. It remains to show that it is surjective. Let $M$ be a maximal ideal of $R$. Suppose $M\neq\mathfrak{m}_{x}\cap R$ for all $x\in X$. Thus for each $x\in X$ there exists some $A_{x}\in R$ such that $A_{x}\in M$ and $A^{c}_{x}\in\mathfrak{m}_{x}\cap R$. It follows that $x\in A_{x}$ for all $x\in X$. But $X$ is quasi-compact, hence we may write $X=\bigcup\limits_{i=1}^{n}A_{x_{i}}$. This yields that $1=X\in M$, see Remark \ref{Remark ii subring new}. But this is a contradiction, since $M$ is a proper ideal of $R$. Hence, $M=\mathfrak{m}_{x}\cap R$ for some $x\in X$. $\Box$

Note that in the proof of Theorem \ref{Theorem V}, we only used the ``quasi-compactness'' assumption of $\pi_{0}(X)$. So we have:

\begin{corollary}\label{Corollary I} If $X$ is a compact space, then $\pi_{0}(X)$ is compact totally disconnected. In particular, it is Hausdorff.
\end{corollary}

{\bf Proof.} It is an immediate consequence of Theorem \ref{Theorem V}. $\Box$

The following well known result (cf. \cite[Proposition 2.15]{Sonia-Sabogal}) is also an immediate consequence of Theorem \ref{Theorem V}. Of course the further relation stated in \cite[Proposition 2.15]{Sonia-Sabogal} is novel.

\begin{corollary}\label{Theorem DIII} Every compact totally disconnected space $X$ is homeomorphic to the prime spectrum of the Boolean ring $\Clop(X)$.
\end{corollary}

{\bf Proof.} It follows from Theorem \ref{Theorem V}. $\Box$

\begin{corollary}\label{Coro iv 2020} If $X$ is a compact totally disconnected space, then the set $\Clop(X)$ forms a base for the opens of $X$.
\end{corollary}

{\bf Proof.} This is deduced from Corollary \ref{Theorem DIII}. Here is another proof. Let $U$ be an open of $X$ and $x\in U$. Then $\{x\}$ is the intersection of all clopens of $X$ containing $x$. Thus there exists a finite number $A_{1},\ldots, A_{n}\in\Clop(X)$ such that $x\in A=\bigcap\limits_{k=1}^{n}A_{k}\subseteq U$, because $U^{c}$ is quasi-compact. Thus $A$ is a clopen of $X$. Hence, $\Clop(X)$ is a base for the opens of $X$. $\Box$

\begin{theorem}\label{Theorem DII} Let $R$ and $R'$ be two Boolean rings. Then the map: $$\Mor_{\Ring}(R,R')\rightarrow\Mor_{\Top}\big(\Spec(R'),\Spec(R)\big)$$ given by $\phi\rightsquigarrow\Spec(\phi)$ is a one-to-one correspondence.
\end{theorem}

{\bf Proof.} Let $\phi,\psi:R\rightarrow R'$ be two ring maps such that $\Spec(\phi)=\Spec(\psi)$. To prove $\phi=\psi$ it suffices to show that for each $e\in R$, then $\phi(e)-\psi(e)$ is a member of the Jacobson radical of $R'$, because it is the zero ideal. Suppose there exists a maximal ideal $\mathfrak{m}$ of $R'$ such that $\phi(e)-\psi(e)\notin\mathfrak{m}$. It follows that $1-\big(\phi(e)-\psi(e)\big)\in\mathfrak{m}$. This yields that $\phi(e)\psi(e)\in\mathfrak{m}$. We may assume that $\phi(e)\in\mathfrak{m}$. This implies that $\psi(e)\notin\mathfrak{m}$. But we have $e\in\phi^{-1}(\mathfrak{m})=\psi^{-1}(\mathfrak{m})$. This is a contradiction and we win. To see the surjectivity, let $h:\Spec(R')\rightarrow\Spec(R)$ be a continuous function. If $e\in R$ then $\D(e)$ is a clopen of $\Spec(R)$ and so $h^{-1}\big(\D(e)\big)$ is a clopen of $\Spec(R')$. Hence there exists a unique $e'\in R'$ such that $h^{-1}\big(\D(e)\big)=\D(e')$. Then consider the function $\phi:R\rightarrow R'$ defined by $e\rightsquigarrow e'$. This map is clearly multiplicative and preserves the unit. To prove its additivity it suffices to show that $h^{-1}\big(\D(e_{1}+e_{2})\big)=\D(e_{1}'+e_{2}')$. If $\mathfrak{p}\in h^{-1}\big(\D(e_{1}+e_{2})\big)$ then
$1-e_{1}-e_{2}\in h(\mathfrak{p})$. Suppose $e_{1}'+e_{2}'\in\mathfrak{p}$. We have $e_{1}e_{2}\in h(\mathfrak{p})$. Thus we may assume that, say, $e_{1}\in h(\mathfrak{p})$. It follows that $e_{2}\notin h(\mathfrak{p})$. Thus $\mathfrak{p}\in h^{-1}\big(\D(e_{2})\big)=\D(e_{2}')$. But $(1-e_{1}')e_{2}'\in\mathfrak{p}$. Therefore $1-e_{1}'\in\mathfrak{p}$. This yields that $\mathfrak{p}\in\D(e_{1}')$. But this is a contradiction. Hence, $h^{-1}\big(\D(e_{1}+e_{2})\big)\subseteq\D(e_{1}'+e_{2}')$. The reverse inclusion is also proved by a similar argument. Thus $\phi$ is a morphism of rings. It remains to show that $\Spec(\phi)=h$. If $\mathfrak{p}$ is a prime ideal of $R'$ then $\phi^{-1}(\mathfrak{p})\subseteq h(\mathfrak{p})$ and so $\phi^{-1}(\mathfrak{p})=h(\mathfrak{p})$, since $\phi^{-1}(\mathfrak{p})$ is a maximal ideal of $R$. $\Box$

Note that the above theorem does not hold in general. For example, if $K$ is a field then every function $\Spec(K)\rightarrow\Spec(\mathbb{Z})$ is continuous, but there exists only one morphism of rings $\mathbb{Z}\rightarrow K$. \\
It is well known and folklore (due to Grothendieck) that the category of (commutative) rings by the contravarinat functor $R\rightsquigarrow\big(\Spec(R),\mathscr{O}_{\Spec(R)}\big)$ and $\phi\rightsquigarrow\big(\Spec(\phi),\phi^{\sharp})$ is anti-equivalent (dual) to the category of affine schemes. Under this duality, the category of zero dimensional rings is anti-equivalent to the category of compact schemes, see \cite[Corollary 3.5]{Tarizadeh-Aghajani}. Moreover, the category of absolutely flat (von Neumann regular) rings is anti-equivalent to the category of compact reduced schemes. All of the above dualities are in the geometric level and do not hold in the topological level. In fact, in the topological level we have Theorems \ref{Corollary V} and \ref{Thm anti-equivalence}.

\begin{theorem}$($Stone duality$)$\label{Corollary V} The category of Boolean rings by the assignments $R\rightsquigarrow\Spec(R)$ and $\phi\rightsquigarrow\Spec(\phi)$ is anti-equivalent to the category of compact totally disconnected spaces.
\end{theorem}

{\bf Proof.} This functor by Theorem \ref{Theorem DII}, is fully-faithful and by Corollary \ref{Theorem DIII}, is essentially surjective. $\Box$

\begin{remark} Here a second proof is given for Theorem \ref{Corollary V} without using Theorem \ref{Theorem DII}. The isomorphism of Corollary \ref{Corollary III}, is functorial. Hence, the identity functor of the category of Boolean rings is isomorphic to $\Clop\circ\Spec$. The isomorphism of Corollary \ref{Theorem DIII}, is also functorial. Therefore the identity functor of the category of compact totally disconnected spaces is isomorphic to $\Spec\circ\Clop$.
\end{remark}

\section{Complete Boolean rings and Stone type duality}

Let $R$ be a Boolean ring. Then we may define a partial order $\leqslant$ over $R$ as $e\leqslant e'$ if $e=ee'$, or equivalently, $\D(e)\subseteq\D(e')$. (Note that this relation $\leqslant$ is a partial order over a ring if and only if it is a Boolean ring). Remember that a Boolean ring $R$ is called complete if it is join-complete with respect to $\leqslant$. It is easy to see that $R$ is join-complete if and only if $R$ is meet-complete with respect to $\leqslant$. In fact, if $(e_{k})\subseteq R$ then $\bigwedge\limits_{k}e_{k}=e$ if and only if $1-e=\bigvee\limits_{k}(1-e_{k})$. \\
Recall that a topological space $X$ is called extremally disconnected if the closure of every open of $X$ is again an open of $X$. It is easy to see that this notion is auto-dual. Namely, a space $X$ is extremally disconnected if and only if the interior of every closed of $X$ is again a closed of $X$. Every irreducible space is extremally disconnected, since in such space every non-empty open is dense. \\
Then we give a new ring-theoretic proof to the following result (which is already known for Boolean algebras).

\begin{lemma}\label{Thm complete Boolean ring} Let $R$ be a Boolean ring. Then $R$ is complete if and only if $\Spec(R)$ is extremally disconnected.
\end{lemma}

{\bf Proof.} Assume $R$ is complete. If $U$ is an open of $\Spec(R)$ then there exists $(e_{k})\subseteq R$ such that $U=\bigcup\limits_{k}\D(e_{k})$. We show that $\overline{U}=\D(e)$ where $e=\bigvee\limits_{k}e_{k}$. Take $\mathfrak{p}\in\overline{U}$. If $e\in\mathfrak{p}$ then $\mathfrak{p}\in\D(1-e)$ and so $\D(1-e)\cap U\neq\emptyset$. Thus $\D\big(e_{k}(1-e)\big)\neq\emptyset$ for some $k$. But this is a contradiction, since $e_{k}(1-e)=0$. Hence, $\overline{U}\subseteq\D(e)$. To see the reverse inclusion, take $\mathfrak{p}\in\D(e)$. It suffices to show that if $\mathfrak{p}\in\D(e')$ for some $e'\in R$, then $\D(e')$ meets $U$. If $U\cap\D(e')=\emptyset$ then $U\subseteq\D(e')^{c}=\V(e')=\D(1-e')$. This shows that $e_{k}\leqslant 1-e'$ for all $k$. Thus $e=e(1-e')\in\mathfrak{p}$ which is a contradiction. Conversely, if $(e_{k})\subseteq R$ then $\overline{U}$ is an open of $\Spec(R)$ where $U=\bigcup\limits_{k}\D(e_{k})$. But $\overline{U}$ is quasi-compact, thus we may write $\overline{U}=\bigcup\limits_{i=1}^{n}\D(e'_{i})=
\bigcup\limits_{i=1}^{n}\V(1-e'_{i})=\V(e')=\D(1-e')$ where $e'=\prod\limits_{i=1}^{n}(1-e'_{i})$. Thus $e_{k}\leqslant 1-e'$ for all $k$. We show that $1-e'=\bigvee\limits_{k}e_{k}$. Suppose for each $k$, $e_{k}\leqslant e''$ for some $e''\in R$. Then $U\subseteq\D(e'')=\V(1-e'')$ and so $\D(1-e')=\overline{U}\subseteq\D(e'')$. $\Box$

\begin{corollary}\label{Theorem extremally} For any set $X$, then the space $\Spec\mathcal{P}(X)$ is extremally disconnected.
\end{corollary}

{\bf Proof.} It is an immediate consequence of Lemma \ref{Thm complete Boolean ring}. $\Box$

By the category of complete Boolean rings we mean the category whose objects are the complete Boolean rings and  whose morphisms are the ring maps which preserve suprema, i.e., morphisms of rings $\phi:R\rightarrow R'$ such that $\phi(\bigvee\limits_{k}e_{k})=\bigvee\limits_{k}\phi(e_{k})$. Clearly a morphism of rings between complete Boolean rings preserves suprema if and only if it preserves infima. Also note that if $R$ and $R'$ are complete Boolean rings such that $R$ is a subring of $R'$ and $(e_{k})$ is a subset of $R$, then $\bigwedge\limits_{k}e_{k}\leqslant\inf\limits_{k}e_{k}$ where
$\bigwedge\limits_{k}e_{k}$ and $\inf\limits_{k}e_{k}$ are the infimums of $(e_{k})$ in $R$ and $R'$, respectively.
In the category of compact extremally disconnected spaces morphisms are the open continuous maps. The following technical result paves the way in order to establish the duality (anti-equivalence) between these two categories.

\begin{theorem}\label{Theorem Guram} Let $\phi:R\rightarrow R'$ be a morphism of rings between complete Boolean rings. Then $\phi$ preserves suprema if and only if the induced map $\phi^{\ast}=\Spec(\phi):\Spec(R')\rightarrow\Spec(R)$ is an open map.
\end{theorem}

{\bf Proof.} Assume $\phi$ preserves suprema. It will be enough to show that for each $U\in\Clop(X_{R'})$ then $\phi^{\ast}(U)$ is an open of $X_{R}=\Spec(R)$ where $X_{R'}=\Spec(R')$. By the hypothesis and Corollary \ref{Corollary III}, $\Clop(\phi^{\ast}):\Clop(X_{R})\rightarrow\Clop(X_{R'})$ preserves infima. Thus there exists an order preserving map $\theta:\Clop(X_{R'})\rightarrow\Clop(X_{R})$ such that the pair $\big(\theta, \Clop(\phi^{\ast})\big)$ is a Galois connection. To prove the assertion it suffices to show that $\phi^{\ast}(U)=\theta(U)$. Clearly $\theta(U)\subseteq\theta(U)$ and so $U\subseteq(\phi^{\ast})^{-1}\big(\theta(U)\big)$, or equivalently,  $\phi^{\ast}(U)\subseteq\theta(U)$. To see the reverse inclusion we act as follows. Clearly $\phi^{\ast}(U)$ is a closed subset of $X_{R}$, since $\phi^{\ast}$ is a closed map. Using Corollary \ref{Coro iv 2020}, then we obtain that $\phi^{\ast}(U)=\bigcap\limits_{V\in\mathcal{S}}V$ where
$\mathcal{S}=\{V\in\Clop(X_{R}): \phi^{\ast}(U)\subseteq V\}$. On the other hand, we have
$\theta(U)=\bigwedge\limits_{V\in\mathcal{S}}V$. Therefore $\theta(U)\subseteq\phi^{\ast}(U)$, because $\Clop(X_{R})$ is a subring of $\mathcal{P}(X_{R})$. Conversely, if $\phi^{\ast}$ is an open map then we obtain a map $\psi:\Clop(X_{R'})\rightarrow\Clop(X_{R})$ given by $U\rightsquigarrow\phi^{\ast}(U)$. Clearly it is an order preserving map and the pair $\big(\psi,\Clop(\phi^{\ast})\big)$ is a Galois connection.
Hence, $\Clop(\phi^{\ast})$ and so $\phi$ preserve infima. $\Box$

\begin{corollary}\label{Coro vi power set} For any function $f:X\rightarrow Y$, the induced map $\Spec\big(\mathcal{P}(f)\big):
\Spec\mathcal{P}(X)\rightarrow
\Spec\mathcal{P}(Y)$ is an open map.
\end{corollary}

{\bf Proof.} It follows from Theorem \ref{Theorem Guram}. $\Box$ \\

We give an alternative proof to the following Stone type duality.

\begin{theorem}\label{Thm anti-equivalence}\cite[Corollary 6.10(2)]{Bezhanishvili} The category of complete Boolean rings by the assignments $R\rightsquigarrow\Spec(R)$ and $\phi\rightsquigarrow\Spec(\phi)$ is anti-equivalent to the category of compact extremally disconnected spaces.
\end{theorem}

{\bf Proof.} If $R$ is a complete Boolean ring then by Lemma \ref{Thm complete Boolean ring}, $\Spec(R)$ is extremally disconnected.
If moreover $\phi$ is a morphism of complete Boolean rings then by Theorem \ref{Theorem Guram}, $\Spec(\phi)$ is an open map. Hence, the above assignments form a contravariant functor. The above functor is fully-faithful and essentially surjective, see Theorems \ref{Corollary V}, \ref{Theorem Guram} and taking into account that every Hausdorff extremally disconnected space is clearly totally disconnected. $\Box$

\begin{remark}\label{Remark i 2020} It is well known that in the category of Hausdorff spaces, epimorphisms are precisely the continuous maps with dense images. By a similar argument, it is shown that in the category of compact spaces, epimorphisms are precisely the surjective continuous maps. It is also well known that a compact space $X$ is extremally disconnected if and only if every continuous and surjective map $f:Y\rightarrow X$ with $Y$ compact admits a continuous section, that is a continuous map $g:X\rightarrow Y$ such that $f\circ g$ is the identity, for its proof see \cite[Tag 08YN]{Johan} or \cite[Theorem 2.5]{Gleason}, also see \cite{Strauss}. Therefore in the category of
compact spaces, the projective objects are precisely the compact extremally disconnected spaces.
\end{remark}

The following result provides a criterion for compact extremally disconnected spaces in terms of the retractions. For the definition of retraction space see e.g. \cite[\S2]{Tarizadeh-Aghajani}.

\begin{theorem}\label{Theorem VI nt} Let $X$ be a compact space. Then $X$ is extremally disconnected if and only if $X$ is a retraction of $\Spec\mathcal{P}(X)$.
\end{theorem}

{\bf Proof.} Assume $X$ is extremally disconnected. By the universal property of the Stone-\v{C}ech compactification, there exists a (unique) continuous map $f:\Spec\mathcal{P}(X)\rightarrow X$ such that $f(\mathfrak{m}_{x})=x$ for all $x\in X$. Thus by Remark \ref{Remark i 2020}, there exists a continuous map $g:X\rightarrow\Spec\mathcal{P}(X)$ such that $f\circ g$ is the identity. Conversely, by Corollary \ref{Theorem extremally}, $\Spec\mathcal{P}(X)$ is extremally disconnected. It is well known that each retraction of a compact extremally disconnected space is extremally disconnected, see \cite[Proposition 24.2.7]{Semandeni}. Hence, $X$ is extremally disconnected. $\Box$

We give a new proof to the following well known result.

\begin{corollary}\label{Coro vii injective} In the category of Boolean rings, the injective objects are precisely the complete Boolean rings.
\end{corollary}

{\bf Proof.} It will be enough to show that in the category of Boolean rings, monomorphisms are precisely the injective ring maps. Because then the assertion is deduced from Remark \ref{Remark i 2020} and Theorem \ref{Theorem VI nt}. So let $\phi:R\rightarrow R'$ be a monomorphism in the category of Boolean rings. Suppose $\phi(e)=\phi(e')$ for some $e,e'\in R$. Assume $e\neq e'$. Setting $X=\{e,e'\}$. Then consider the function $f:\mathcal{P}(X)\rightarrow R$ defined by $f(\emptyset)=0$, $f(X)=1$, $f(\{e\})=e$ and $f(\{e'\})=1-e$. Similarly we define $g:\mathcal{P}(X)\rightarrow R$ as $g(\emptyset)=0$, $g(X)=1$, $g(\{e\})=e'$ and $g(\{e'\})=1-e'$. Clearly $f$ and $g$ are morphisms of rings and $\phi\circ f=\phi\circ g$. It follows that $f=g$ which is a contradiction. $\Box$

We say that a morphism of rings $\phi:R\rightarrow R'$ admits a right inverse if there exists a morphism of rings $\psi:R'\rightarrow R$ such that $\phi\circ\psi$ is the identity. The left inverse notion is defined dually. The following result provides a new criterion for the completeness of Boolean rings.

\begin{theorem}\label{Coro viii complete} Let $R$ be a ring and $X=\Spec(R)$. Then the following statements are equivalent. \\
$\mathbf{(i)}$ $\mathcal{B}(R)$ is complete. \\
$\mathbf{(ii)}$ The canonical ring map $\psi:\mathcal{B}(R)\rightarrow\mathcal{P}(X)$ given by $e\rightsquigarrow D(e)$ admits a left inverse.
\end{theorem}

{\bf Proof.} $\mathbf{(i)}\Rightarrow\mathbf{(ii)}:$ By Corollary \ref{Corollary ix injective b}, $\psi$ is injective. Then apply Corollary \ref{Coro vii injective}. \\
$\mathbf{(ii)}\Rightarrow\mathbf{(i)}:$ By Lemma \ref{Thm complete Boolean ring}, it suffices to show that $Y=\Spec\mathcal{B}(R)$ is extremally disconnected. Then by Theorem \ref{Theorem VI nt}, it will be enough to show that $Y$ is a retraction of $\Spec\mathcal{P}(Y)$. By the hypothesis, there exists a morphism of rings $\phi:\mathcal{P}(X)\rightarrow\mathcal{B}(R)$ such that $\phi\circ\psi$ is the identity map. Setting $\theta:=\phi\circ\mathcal{P}(\mu^{-1}\circ\pi)$ where the canonical projection map $\pi:X\rightarrow \Sp(R)$ is defined as $\mathfrak{p}\rightsquigarrow\mathfrak{p}^{\ast}=(e\in\mathfrak{p}: e=e^{2})$ and for $\mu:Y\rightarrow\Sp(R)$ see Theorem \ref{Theorem PSp}. One can observe that $\psi=\mathcal{P}(\mu^{-1}\circ\pi)\circ\psi'$ where $\psi':\mathcal{B}(R)\rightarrow\mathcal{P}(Y)$ is the canonical ring map which is given by $e\rightsquigarrow D(e)$. Clearly $\Spec(\psi')\circ\Spec(\theta)$ is the identity map. Hence,
$Y$ is a retraction of $\Spec\mathcal{P}(Y)$. $\Box$

\begin{corollary} Let $R$ be a Boolean ring and $X=\Spec(R)$. Then the following statements are equivalent. \\
$\mathbf{(i)}$ $R$ is complete. \\
$\mathbf{(ii)}$ There exists a morphism of rings $\phi:\mathcal{P}(X)\rightarrow R$ which admits a right inverse.
\end{corollary}

{\bf Proof.} $\mathbf{(i)}\Rightarrow\mathbf{(ii)}:$ See Theorem \ref{Coro viii complete}(ii). \\
$\mathbf{(ii)}\Rightarrow\mathbf{(i)}:$ By the hypothesis, $X$ is a retraction of $\Spec\mathcal{P}(X)$. Thus by Theorem \ref{Theorem VI nt}, $X$ is extremally disconnected. So by Lemma \ref{Thm complete Boolean ring}, $R$ is complete. $\Box$

\section{Power set ring and fixed-point theorems}

\begin{proposition}\label{Remark III} Let $M$ be a maximal ideal of a Boolean ring $R$. Then the map $\phi_{M}:R\rightarrow\mathbb{Z}_{2}$ which sends each $e\in R$ into $0$ or $1$, according as $e\in M$ or $e\notin M$, is a morphism of rings.
\end{proposition}

{\bf Proof.} Clearly $\phi_{M}$ is multiplicative and preserves the identity. To prove the additivity of $\phi_{M}$, it suffices to show that if $e,e'\notin M$ then $e+e'\in M$. If $e+e'\notin M$ then $e+ee'=e(e+e')\notin M$, similarly we get that $e'+ee'=e'(e+e')\notin M$. It follows that $0=(e+ee')(e'+ee')\notin M$, which is impossible. $\Box$

\begin{corollary}\label{Proposition II} If $R$ is a Boolean ring, then we have the canonical bijection $\Spec(R)\simeq\Mor_{\Ring}(R,\mathbb{Z}_{2})$.
\end{corollary}

{\bf Proof.} The map $\Spec(R)\rightarrow\Mor_{\Ring}(R,\mathbb{Z}_{2})$ given by $M\rightsquigarrow\phi_{M}$ is clearly injective, for $\phi_{M}$ see Proposition \ref{Remark III}. It is also surjective, because if $\theta:R\rightarrow\mathbb{Z}_{2}$ is a morphism of rings then $\theta=\phi_{M}$ where $M=\theta^{-1}(0)$. Here is another proof, the assertion is a special case of Theorem \ref{Theorem DII}. $\Box$

Let $X$ be a set with the cardinality $\kappa$. Remember that $\Spec\mathcal{P}(X)$ is the Stone-\v{C}ech compactification of the discrete space $X$. Using this and Corollary \ref{Proposition II}, then the number of all ring maps $\mathcal{P}(X)\rightarrow\mathbb{Z}_{2}$ is either $\kappa$ or $2^{2^{\kappa}}$, according as $X$ is finite or infinite.

\begin{remark} The power set functor is not full, see Theorem \ref{Corollary II}. Let $\phi:\mathcal{P}(X)\rightarrow\mathcal{P}(Y)$ be a morphism of rings. Then it is natural to ask, under which circumstances there is a function $f:Y\rightarrow X$ such that $\phi=\mathcal{P}(f)$? The answer is affirmative precisely when $X$ is finite, see Theorem \ref{Corollary II}.
\end{remark}

Let $X$ be an infinite set and let $\Fin(X)$ be the set of all finite subsets of $X$. Then $\Fin(X)$ is an ideal of $\mathcal{P}(X)$ which is not a finitely generated ideal. Because if it is a finitely generated ideal then it is principal, say $\Fin(X)=(A)$. Clearly $X\setminus A$ is non-empty, so choose some $x$ in it. Then $\{x\}\in\Fin(X)$ and so $\{x\}=BA=B\cap A$ for some $B\in\mathcal{P}(X)$, but this is a contradiction. Hence, $\Fin(X)$ is not a finitely generated ideal. Indeed, $\Fin(X)$ is generated by the single-point subsets of $X$, since if $A\in\Fin(X)$ then we may write $A=\sum\limits_{x\in A}\{x\}$.

\begin{theorem}\label{Corollary II} For a set $X$ the following statements are equivalent. \\
$\mathbf{(i)}$ $X$ is a finite set. \\
$\mathbf{(ii)}$ The maximal ideals of $\mathcal{P}(X)$ are precisely of the form $\mathfrak{m}_{x}$. \\
$\mathbf{(iii)}$ For each nonempty set $Y$, then every morphism of rings $\mathcal{P}(X)\rightarrow\mathcal{P}(Y)$ is of the form $\mathcal{P}(f)$ for some function $f:Y\rightarrow X$.
\end{theorem}

{\bf Proof.} $\mathbf{(i)}\Rightarrow\mathbf{(ii)}:$ Let $M$ be a maximal ideal of $\mathcal{P}(X)$. We have $\bigcap\limits_{x\in X}\mathfrak{m}_{x}=0$. But $X$ is finite, thus $\mathfrak{m}_{x}\subseteq M$
for some $x\in X$. So $\mathfrak{m}_{x}=M$. \\
$\mathbf{(ii)}\Rightarrow\mathbf{(iii)}:$ Let $\phi:\mathcal{P}(X)\rightarrow\mathcal{P}(Y)$ be a morphism of rings. If $y\in Y$ there exists a unique point $x\in X$ such that $\phi^{-1}(\mathfrak{m}_{y})=\mathfrak{m}_{x}$. Then the map $f:Y\rightarrow X$ defined by $y\rightsquigarrow x$ is the desired function. \\
$\mathbf{(iii)}\Rightarrow\mathbf{(i)}:$ If $X$ is an infinite set, then $\Fin(X)$ is a proper ideal of $\mathcal{P}(X)$. Thus there exists a maximal ideal $M$ of $\mathcal{P}(X)$ such that $\Fin(X)\subseteq M$. Then the map $\phi:\mathcal{P}(X)\rightarrow\mathcal{P}(Y)$ which sends each $A\in\mathcal{P}(X)$ into $0$ or $1$, according as $A\in M$ or $A\notin M$, is a morphism of rings, see Proposition \ref{Remark III}. Thus by the hypothesis, there exists a function $f:Y\rightarrow X$ such that $\phi=\mathcal{P}(f)$.
But the image of $f$ is nonempty, hence we may choose some point $x$ in it. We have $\{x\}\in M$ and so $\phi(\{x\})=0$. But $\phi(\{x\})=f^{-1}(\{x\})\neq\emptyset$. This is a contradiction. Hence, $X$ is finite. $\Box$

\begin{theorem}\label{Theorem II} Let $M$ be a maximal ideal of $\mathcal{P}(X)$. Then the following statements are equivalent. \\
$\mathbf{(i)}$ $\Fin(X)\subseteq M$. \\
$\mathbf{(ii)}$ $M\neq\mathfrak{m}_{x}$ for all $x\in X$. \\
$\mathbf{(iii)}$ $M$ is not a finitely generated ideal.
\end{theorem}

{\bf Proof.} $\mathbf{(i)}\Rightarrow\mathbf{(ii)}:$ If $M=\mathfrak{m}_{x}=\mathcal{P}(X\setminus\{x\})$
for some $x\in X$ then $\{x\}\notin M$, which is a contradiction. \\
$\mathbf{(ii)}\Rightarrow\mathbf{(iii)}:$ If $M$ is a finitely generated ideal then it is principal. Thus $M=(A)=\mathcal{P}(A)$ for some $A\in\mathcal{P}(X)$. But $A$ is a proper subset of $X$, hence we may choose some $x\in X\setminus A$. Thus $\{x\}\notin M$ and so $X\setminus\{x\}\in M$. This yields that $M=\mathfrak{m}_{x}$. But this is a contradiction. \\
$\mathbf{(iii)}\Rightarrow\mathbf{(i)}:$ It will be enough to show that $\{x\}\in M$ for all $x\in X$. Suppose there is some $x\in X$ such that $\{x\}\notin M$. It follows that $X\setminus\{x\}\in M$ and so $M=\mathfrak{m}_{x}$ is a principal ideal. But this is a contradiction. $\Box$

Theorem \ref{Theorem II}, in particular, tells us that the maximal ideals of $\mathcal{P}(X)$ are either infinitely generated or the principal ideals of the form $\mathfrak{m}_{x}$. \\
Let $X$ be a set and let $M$ be a maximal ideal of $\mathcal{P}(X)$. Then by Theorem \ref{Theorem II}, $M+\Fin(X)=\mathcal{P}(X)$ if and only if $M=\mathfrak{m}_{x}$ for some $x\in X$. \\
Theorem \ref{Corollary II} also leads us to the following general conclusion.

\begin{corollary} Let $X$ and $Y$ be two sets. Then the image of the induced map $\Mor_{\Set}(X,Y)\rightarrow\Mor_{\Ring}
\big(\mathcal{P}(Y),\mathcal{P}(X)\big)$ given by
$f\rightsquigarrow\mathcal{P}(f)$ is consisting of all morphisms of rings $\phi:\mathcal{P}(Y)\rightarrow\mathcal{P}(X)$ such that $\Fin(Y)+\phi^{-1}(\mathfrak{m}_{x})=\mathcal{P}(Y)$ for all $x\in X$.
\end{corollary}

{\bf Proof.} First we need to show that if $f:X\rightarrow Y$ is a function then the induced ring map $\psi=\mathcal{P}(f)$ actually satisfies in the above condition. Suppose there is some $x\in X$ such that $\Fin(Y)\subseteq\psi^{-1}(\mathfrak{m}_{x})$. This in particular yields that $f^{-1}(\{y\})\subseteq X\setminus\{x\}$ which is impossible where $y=f(x)$. Conversely, let $\phi:\mathcal{P}(Y)\rightarrow\mathcal{P}(X)$ be any ring map which satisfies in the above condition. For each $x\in X$ then by Theorem \ref{Theorem II}, there exists a unique point $y\in Y$ such that $\phi^{-1}(\mathfrak{m}_{x})=\mathfrak{m}_{y}$. Now it is easy to see that $\phi=\mathcal{P}(f)$ where the function $f:X\rightarrow Y$ is defined as $x\rightsquigarrow y$. $\Box$

\begin{theorem}\label{Theorem III} If $X$ is an infinite set, then every maximal ideal of the quotient ring $\mathcal{P}(X)/\Fin(X)$ is infinitely generated.
\end{theorem}

{\bf Proof.} Let $M$ be a maximal ideal of $\mathcal{P}(X)$ such that $\Fin(X)\subseteq M$. If $M/\Fin(X)$ is finitely generated then it is a principal ideal. Thus there exists some $A\in\mathcal{P}(X)$ such that  $M=\Fin(X)+\mathcal{P}(A)$. But $A^{c}=X\setminus A$ is an infinite set, because if it is finite then $1=A^{c}+A\in M$, which is impossible. Thus we may choose two infinite and disjoint subsets $B$ and $C$ in the infinite set $A^{c}$. Therefore $BC=0\in M$. But none of them is in $M$, because if for example $B\in M$ then there exist a finite subset $F\subseteq X$ and some $A'\in\mathcal{P}(A)$ such that $B=F+A'$, but $B\cap A'=\emptyset$ and so $B=B^{2}=BF+BA'=B\cap F+B\cap A'=B\cap F\subseteq F$ which is impossible, since $B$ is infinite. $\Box$

If $X$ is an infinite set then by Theorem \ref{Theorem III}, the quotient ring $\mathcal{P}(X)/\Fin(X)$ is not a Noetherian ring. In particular, $\Fin(X)$ is not a maximal ideal of $\mathcal{P}(X)$. \\
In the following results we have given fresh proofs of the known facts, see e.g. \cite[\S3.4]{Hindman-Strauss}. These proofs are simpler than the original arguments, because we use the elementary and standard language of commutative algebra both in their formulations and the methods used in their proofs. \\
The following result provides a general criterion guaranteeing that, if it is satisfied, then a function yields a fixed-point.

\begin{theorem}\label{Theorem DVI} A function $f:X\rightarrow X$ has a fixed-point if and only if the induced map
$\Spec\big(\mathcal{P}(f)\big):
\Spec\mathcal{P}(X)\rightarrow\Spec\mathcal{P}(X)$ has a fixed-point.
\end{theorem}

{\bf Proof.} Setting $\phi=\Spec\big(\mathcal{P}(f)\big)$. If $f(x)=x$ for some $x\in X$ then $\phi(\mathfrak{m}_{x})=\mathfrak{m}_{x}$. Conversely, assume that $\phi(M)=M$ for some $M\in\Spec\mathcal{P}(X)$. If $f:X\rightarrow X$ does not have any fixed-point, then it is well known that there exists a function $g:X\rightarrow\{1,2,3\}$ such that $g(x)\neq g\big(f(x)\big)$ for all $x\in X$. (The existence of such function is the key point of our proof, so we sketch its proof. First, Zorn's lemma gives us a maximal function $g:A\rightarrow\{1,2,3\}$ satisfying in the above condition where $A\in\mathcal{P}(X)$ with $f(A)\subseteq A$. Then we show that $A=X$. If $x\in X\setminus A$, then by the induction we may define a function $h:A\cup\{x,f(x), f^{2}(x),\ldots\}\rightarrow\{1,2,3\}$ extending $g$ and still satisfying in the above condition as follows, for each natural number $n\geqslant0$ we choose $h\big(f^{n}(x)\big)$ to be a value in $\{1,2,3\}$ in a way that $h\big(f^{n}(x)\big)\neq h\big(f^{n-1}(x)\big)$ and $h\big(f^{n}(x)\big)\neq h\big(f^{n+1}(x)\big)$ if these values have already been defined. Finally, we reach to a contradiction with the maximality of $g$, hence $A=X$). Now setting $C_{k}=g^{-1}(\{k\})$. Then it is obvious that $C_{k}\cap f(C_{k})=\emptyset$ for $k=1,2,3$. But $1\notin M$ and so $C_{k}\notin M$ for some $k$. It follows that $f(C_{k})\in M$. We have then $C_{k}\subseteq f^{-1}\big(f(C_{k})\big)\in M$. Thus $C_{k}\in M$.
But this is a contradiction. Hence, $f$ has a fixed-point.  $\Box$

As an application of Theorem \ref{Theorem DVI}, the following result is obtained.

\begin{theorem}\label{Theorem fixed point} Let $f:X\rightarrow X$ be a function, $\phi=\Spec\big(\mathcal{P}(f)\big)$ and $M\in\Spec\mathcal{P}(X)$. Then $\phi(M)=M$ if and only if $\{x\in X: f(x)\neq x\}\in M$.
\end{theorem}

{\bf Proof.} Let $A=\{x\in X: f(x)\neq x\}$. If $A\in M$ then to prove the assertion it suffices to show that $M\subseteq\phi(M)$, because $M$ is a maximal ideal. Suppose there exists some $B\in M$ such that $B\notin\phi(M)$. It follows that $f^{-1}(B)\notin M$. Clearly $ f^{-1}(B)\cap A^{c}\subseteq B$. Therefore $f^{-1}(B)\cap A^{c}\in M$. But this is a contradiction. Conversely, assume that $\phi(M)=M$. If $A\notin M$ then $A$ is non-empty, hence we may choose a point $p$ in it. Then consider the function $g:X\rightarrow X$ which sends each $x\in X$ into $f(x)$ or $p$, according as $x\in A$ or $x\in A^{c}$. Therefore the induced map $\psi=\Spec\big(\mathcal{P}(g)\big)$ and $\phi$ agree on $\overline{\eta(A)}$ where $\eta:X\rightarrow\Spec\mathcal{P}(X)$ is the canonical map which is defined as $\eta(x)=\mathfrak{m}_{x}$. But $M\in\overline{\eta(A)}$, because if $M\in\D(B)$ for some $B\in\mathcal{P}(X)$ then $A\cap B\notin M$ and hence $A\cap B$ is nonempty, thus we may choose some $x\in A\cap B$, this yields that $\mathfrak{m}_{x}\in\D(B)\cap\eta(A)$. So $\psi(M)=M$. Thus by Theorem \ref{Theorem DVI}, $g$ has a fixed-point. But this is a contradiction, since $g$ has no fixed-point. $\Box$ \\

\textbf{Acknowledgements.} The authors would like to give heartfelt thanks to the professors Pierre Deligne and Guram Bezhanishvili for the useful discussions and for their very valuable suggestions and comments which greatly improved the paper. We would also like to give sincere thanks to the referee for very careful reading of the paper. \\

\end{document}